\input amstex
\input amsppt.sty

\magnification1200

\hsize13cm
\vsize19cm

\TagsOnRight

\catcode`\@=11
\font\tenln    = line10
\font\tenlnw   = linew10

\newskip\Einheit \Einheit=0.5cm
\newcount\xcoord \newcount\ycoord
\newdimen\xdim \newdimen\ydim \newdimen\PfadD@cke \newdimen\Pfadd@cke

%%%%%%%%%%%%%%%%%%%%%%%%%%%%%%%%%%%%%%%%%%%%%%%%%
%LaTeX counters, dimensions, variables for lines%
%%%%%%%%%%%%%%%%%%%%%%%%%%%%%%%%%%%%%%%%%%%%%%%%%
\newcount\@tempcnta
\newcount\@tempcntb

\newdimen\@tempdima
\newdimen\@tempdimb

\newdimen\@wholewidth
\newdimen\@halfwidth

\newcount\@xarg
\newcount\@yarg
\newcount\@yyarg
\newbox\@linechar
\newbox\@tempboxa
\newdimen\@linelen
\newdimen\@clnwd
\newdimen\@clnht

\newif\if@negarg

\def\@whilenoop#1{}
\def\@whiledim#1\do #2{\ifdim #1\relax#2\@iwhiledim{#1\relax#2}\fi}
\def\@iwhiledim#1{\ifdim #1\let\@nextwhile=\@iwhiledim
        \else\let\@nextwhile=\@whilenoop\fi\@nextwhile{#1}}

\def\@whileswnoop#1\fi{}
\def\@whilesw#1\fi#2{#1#2\@iwhilesw{#1#2}\fi\fi}
\def\@iwhilesw#1\fi{#1\let\@nextwhile=\@iwhilesw
         \else\let\@nextwhile=\@whileswnoop\fi\@nextwhile{#1}\fi}

\def\thinlines{\let\@linefnt\tenln \let\@circlefnt\tencirc
  \@wholewidth\fontdimen8\tenln \@halfwidth .5\@wholewidth}
\def\thicklines{\let\@linefnt\tenlnw \let\@circlefnt\tencircw
  \@wholewidth\fontdimen8\tenlnw \@halfwidth .5\@wholewidth}
\thinlines
%%%%%%%%%%%%%%%%%%%%%%%%%%%%%%%%%%%%%%%%%%%%%%%%%%%%%%%%%%%

\PfadD@cke1pt \Pfadd@cke0.5pt
\def\PfadDicke#1{\PfadD@cke#1 \divide\PfadD@cke by2 \Pfadd@cke\PfadD@cke \multiply\PfadD@cke by2}
\long\def\LOOP#1\REPEAT{\def\BODY{#1}\ITERATE}
\def\ITERATE{\BODY \let\next\ITERATE \else\let\next\relax\fi \next}
\let\REPEAT=\fi
\def\Punkt{\hbox{\raise-2pt\hbox to0pt{\hss$\ssize\bullet$\hss}}}
\def\DuennPunkt(#1,#2){\unskip
  \raise#2 \Einheit\hbox to0pt{\hskip#1 \Einheit
          \raise-2.5pt\hbox to0pt{\hss$\bullet$\hss}\hss}}
\def\NormalPunkt(#1,#2){\unskip
  \raise#2 \Einheit\hbox to0pt{\hskip#1 \Einheit
          \raise-3pt\hbox to0pt{\hss\twelvepoint$\bullet$\hss}\hss}}
\def\DickPunkt(#1,#2){\unskip
  \raise#2 \Einheit\hbox to0pt{\hskip#1 \Einheit
          \raise-4pt\hbox to0pt{\hss\fourteenpoint$\bullet$\hss}\hss}}
\def\Kreis(#1,#2){\unskip
  \raise#2 \Einheit\hbox to0pt{\hskip#1 \Einheit
          \raise-4pt\hbox to0pt{\hss\fourteenpoint$\circ$\hss}\hss}}

%%%%%%%%%%%%%%%%%%%%%
%LaTeX line macros%
%%%%%%%%%%%%%%%%%%%%%
\def\Line@(#1,#2)#3{\@xarg #1\relax \@yarg #2\relax
\@linelen=#3\Einheit
\ifnum\@xarg =0 \@vline
  \else \ifnum\@yarg =0 \@hline \else \@sline\fi
\fi}

\def\@sline{\ifnum\@xarg< 0 \@negargtrue \@xarg -\@xarg \@yyarg -\@yarg
  \else \@negargfalse \@yyarg \@yarg \fi
\ifnum \@yyarg >0 \@tempcnta\@yyarg \else \@tempcnta -\@yyarg \fi
\ifnum\@tempcnta>6 \@badlinearg\@tempcnta0 \fi
\ifnum\@xarg>6 \@badlinearg\@xarg 1 \fi
\setbox\@linechar\hbox{\@linefnt\@getlinechar(\@xarg,\@yyarg)}%
\ifnum \@yarg >0 \let\@upordown\raise \@clnht\z@
   \else\let\@upordown\lower \@clnht \ht\@linechar\fi
\@clnwd=\wd\@linechar
\if@negarg \hskip -\wd\@linechar \def\@tempa{\hskip -2\wd\@linechar}\else
     \let\@tempa\relax \fi
\@whiledim \@clnwd <\@linelen \do
  {\@upordown\@clnht\copy\@linechar
   \@tempa
   \advance\@clnht \ht\@linechar
   \advance\@clnwd \wd\@linechar}%
\advance\@clnht -\ht\@linechar
\advance\@clnwd -\wd\@linechar
\@tempdima\@linelen\advance\@tempdima -\@clnwd
\@tempdimb\@tempdima\advance\@tempdimb -\wd\@linechar
\if@negarg \hskip -\@tempdimb \else \hskip \@tempdimb \fi
\multiply\@tempdima \@m
\@tempcnta \@tempdima \@tempdima \wd\@linechar \divide\@tempcnta \@tempdima
\@tempdima \ht\@linechar \multiply\@tempdima \@tempcnta
\divide\@tempdima \@m
\advance\@clnht \@tempdima
\ifdim \@linelen <\wd\@linechar
   \hskip \wd\@linechar
  \else\@upordown\@clnht\copy\@linechar\fi}

\def\@hline{\ifnum \@xarg <0 \hskip -\@linelen \fi
\vrule height\Pfadd@cke width \@linelen depth\Pfadd@cke
\ifnum \@xarg <0 \hskip -\@linelen \fi}

\def\@getlinechar(#1,#2){\@tempcnta#1\relax\multiply\@tempcnta 8
\advance\@tempcnta -9 \ifnum #2>0 \advance\@tempcnta #2\relax\else
\advance\@tempcnta -#2\relax\advance\@tempcnta 64 \fi
\char\@tempcnta}

\def\Vektor(#1,#2)#3(#4,#5){\unskip\leavevmode
  \xcoord#4\relax \ycoord#5\relax
      \raise\ycoord \Einheit\hbox to0pt{\hskip\xcoord \Einheit
         \Vector@(#1,#2){#3}\hss}}

\def\Vector@(#1,#2)#3{\@xarg #1\relax \@yarg #2\relax
\@tempcnta \ifnum\@xarg<0 -\@xarg\else\@xarg\fi
\ifnum\@tempcnta<5\relax
\@linelen=#3\Einheit
\ifnum\@xarg =0 \@vvector
  \else \ifnum\@yarg =0 \@hvector \else \@svector\fi
\fi
\else\@badlinearg\fi}

\def\@hvector{\@hline\hbox to 0pt{\@linefnt
\ifnum \@xarg <0 \@getlarrow(1,0)\hss\else
    \hss\@getrarrow(1,0)\fi}}

\def\@vvector{\ifnum \@yarg <0 \@downvector \else \@upvector \fi}

\def\@svector{\@sline
\@tempcnta\@yarg \ifnum\@tempcnta <0 \@tempcnta=-\@tempcnta\fi
\ifnum\@tempcnta <5
  \hskip -\wd\@linechar
  \@upordown\@clnht \hbox{\@linefnt  \if@negarg
  \@getlarrow(\@xarg,\@yyarg) \else \@getrarrow(\@xarg,\@yyarg) \fi}%
\else\@badlinearg\fi}

\def\@upline{\hbox to \z@{\hskip -.5\Pfadd@cke \vrule width \Pfadd@cke
   height \@linelen depth \z@\hss}}

\def\@downline{\hbox to \z@{\hskip -.5\Pfadd@cke \vrule width \Pfadd@cke
   height \z@ depth \@linelen \hss}}

\def\@upvector{\@upline\setbox\@tempboxa\hbox{\@linefnt\char'66}\raise
     \@linelen \hbox to\z@{\lower \ht\@tempboxa\box\@tempboxa\hss}}

\def\@downvector{\@downline\lower \@linelen
      \hbox to \z@{\@linefnt\char'77\hss}}

\def\@getlarrow(#1,#2){\ifnum #2 =\z@ \@tempcnta='33\else
\@tempcnta=#1\relax\multiply\@tempcnta \sixt@@n \advance\@tempcnta
-9 \@tempcntb=#2\relax\multiply\@tempcntb \tw@
\ifnum \@tempcntb >0 \advance\@tempcnta \@tempcntb\relax
\else\advance\@tempcnta -\@tempcntb\advance\@tempcnta 64
\fi\fi\char\@tempcnta}

\def\@getrarrow(#1,#2){\@tempcntb=#2\relax
\ifnum\@tempcntb < 0 \@tempcntb=-\@tempcntb\relax\fi
\ifcase \@tempcntb\relax \@tempcnta='55 \or
\ifnum #1<3 \@tempcnta=#1\relax\multiply\@tempcnta
24 \advance\@tempcnta -6 \else \ifnum #1=3 \@tempcnta=49
\else\@tempcnta=58 \fi\fi\or
\ifnum #1<3 \@tempcnta=#1\relax\multiply\@tempcnta
24 \advance\@tempcnta -3 \else \@tempcnta=51\fi\or
\@tempcnta=#1\relax\multiply\@tempcnta
\sixt@@n \advance\@tempcnta -\tw@ \else
\@tempcnta=#1\relax\multiply\@tempcnta
\sixt@@n \advance\@tempcnta 7 \fi\ifnum #2<0 \advance\@tempcnta 64 \fi
\char\@tempcnta}
%%%%%%%%%%%%%%%%%%%%%%%%%%%%%%%%%%%%%%%%%%%%%%%%%%%%%%%%%%%%%

\def\Diagonale(#1,#2)#3{\unskip\leavevmode
  \xcoord#1\relax \ycoord#2\relax
      \raise\ycoord \Einheit\hbox to0pt{\hskip\xcoord \Einheit
         \Line@(1,1){#3}\hss}}
\def\AntiDiagonale(#1,#2)#3{\unskip\leavevmode
  \xcoord#1\relax \ycoord#2\relax %\advance\xcoord by -0.05\relax
      \raise\ycoord \Einheit\hbox to0pt{\hskip\xcoord \Einheit
         \Line@(1,-1){#3}\hss}}
\def\Pfad(#1,#2),#3\endPfad{\unskip\leavevmode
  \xcoord#1 \ycoord#2 \thicklines\ZeichnePfad#3\endPfad\thinlines}
\def\ZeichnePfad#1{\ifx#1\endPfad\let\next\relax
  \else\let\next\ZeichnePfad
    \ifnum#1=1
      \raise\ycoord \Einheit\hbox to0pt{\hskip\xcoord \Einheit
         \vrule height\Pfadd@cke width1 \Einheit depth\Pfadd@cke\hss}%
      \advance\xcoord by 1
    \else\ifnum#1=2
      \raise\ycoord \Einheit\hbox to0pt{\hskip\xcoord \Einheit
        \hbox{\hskip-\PfadD@cke\vrule height1 \Einheit width\PfadD@cke depth0pt}\hss}%
      \advance\ycoord by 1
    \else\ifnum#1=3
      \raise\ycoord \Einheit\hbox to0pt{\hskip\xcoord \Einheit
         \Line@(1,1){1}\hss}
      \advance\xcoord by 1
      \advance\ycoord by 1
    \else\ifnum#1=4
      \raise\ycoord \Einheit\hbox to0pt{\hskip\xcoord \Einheit
         \Line@(1,-1){1}\hss}
      \advance\xcoord by 1
      \advance\ycoord by -1
    \fi\fi\fi\fi
  \fi\next}
\def\hSSchritt{\leavevmode\raise-.4pt\hbox to0pt{\hss.\hss}\hskip.2\Einheit
  \raise-.4pt\hbox to0pt{\hss.\hss}\hskip.2\Einheit
  \raise-.4pt\hbox to0pt{\hss.\hss}\hskip.2\Einheit
  \raise-.4pt\hbox to0pt{\hss.\hss}\hskip.2\Einheit
  \raise-.4pt\hbox to0pt{\hss.\hss}\hskip.2\Einheit}
\def\vSSchritt{\vbox{\baselineskip.2\Einheit\lineskiplimit0pt
\hbox{.}\hbox{.}\hbox{.}\hbox{.}\hbox{.}}}
\def\DSSchritt{\leavevmode\raise-.4pt\hbox to0pt{%
  \hbox to0pt{\hss.\hss}\hskip.2\Einheit
  \raise.2\Einheit\hbox to0pt{\hss.\hss}\hskip.2\Einheit
  \raise.4\Einheit\hbox to0pt{\hss.\hss}\hskip.2\Einheit
  \raise.6\Einheit\hbox to0pt{\hss.\hss}\hskip.2\Einheit
  \raise.8\Einheit\hbox to0pt{\hss.\hss}\hss}}
\def\dSSchritt{\leavevmode\raise-.4pt\hbox to0pt{%
  \hbox to0pt{\hss.\hss}\hskip.2\Einheit
  \raise-.2\Einheit\hbox to0pt{\hss.\hss}\hskip.2\Einheit
  \raise-.4\Einheit\hbox to0pt{\hss.\hss}\hskip.2\Einheit
  \raise-.6\Einheit\hbox to0pt{\hss.\hss}\hskip.2\Einheit
  \raise-.8\Einheit\hbox to0pt{\hss.\hss}\hss}}
\def\SPfad(#1,#2),#3\endSPfad{\unskip\leavevmode
  \xcoord#1 \ycoord#2 \ZeichneSPfad#3\endSPfad}
\def\ZeichneSPfad#1{\ifx#1\endSPfad\let\next\relax
  \else\let\next\ZeichneSPfad
    \ifnum#1=1
      \raise\ycoord \Einheit\hbox to0pt{\hskip\xcoord \Einheit
         \hSSchritt\hss}%
      \advance\xcoord by 1
    \else\ifnum#1=2
      \raise\ycoord \Einheit\hbox to0pt{\hskip\xcoord \Einheit
        \hbox{\hskip-2pt \vSSchritt}\hss}%
      \advance\ycoord by 1
    \else\ifnum#1=3
      \raise\ycoord \Einheit\hbox to0pt{\hskip\xcoord \Einheit
         \DSSchritt\hss}
      \advance\xcoord by 1
      \advance\ycoord by 1
    \else\ifnum#1=4
      \raise\ycoord \Einheit\hbox to0pt{\hskip\xcoord \Einheit
         \dSSchritt\hss}
      \advance\xcoord by 1
      \advance\ycoord by -1
    \fi\fi\fi\fi
  \fi\next}
\def\Koordinatenachsen(#1,#2){\unskip
 \hbox to0pt{\hskip-.5pt\vrule height#2 \Einheit width.5pt depth1 \Einheit}%
 \hbox to0pt{\hskip-1 \Einheit \xcoord#1 \advance\xcoord by1
    \vrule height0.25pt width\xcoord \Einheit depth0.25pt\hss}}
\def\Koordinatenachsen(#1,#2)(#3,#4){\unskip
 \hbox to0pt{\hskip-.5pt \ycoord-#4 \advance\ycoord by1
    \vrule height#2 \Einheit width.5pt depth\ycoord \Einheit}%
 \hbox to0pt{\hskip-1 \Einheit \hskip#3\Einheit 
    \xcoord#1 \advance\xcoord by1 \advance\xcoord by-#3 
    \vrule height0.25pt width\xcoord \Einheit depth0.25pt\hss}}
\def\Gitter(#1,#2){\unskip \xcoord0 \ycoord0 \leavevmode
  \LOOP\ifnum\ycoord<#2
    \loop\ifnum\xcoord<#1
      \raise\ycoord \Einheit\hbox to0pt{\hskip\xcoord \Einheit\Punkt\hss}%
      \advance\xcoord by1
    \repeat
    \xcoord0
    \advance\ycoord by1
  \REPEAT}
\def\Gitter(#1,#2)(#3,#4){\unskip \xcoord#3 \ycoord#4 \leavevmode
  \LOOP\ifnum\ycoord<#2
    \loop\ifnum\xcoord<#1
      \raise\ycoord \Einheit\hbox to0pt{\hskip\xcoord \Einheit\Punkt\hss}%
      \advance\xcoord by1
    \repeat
    \xcoord#3
    \advance\ycoord by1
  \REPEAT}
\def\Label#1#2(#3,#4){\unskip \xdim#3 \Einheit \ydim#4 \Einheit
  \def\lo{\advance\xdim by-.5 \Einheit \advance\ydim by.5 \Einheit}%
  \def\llo{\advance\xdim by-.25cm \advance\ydim by.5 \Einheit}%
  \def\loo{\advance\xdim by-.5 \Einheit \advance\ydim by.25cm}%
  \def\o{\advance\ydim by.25cm}%
  \def\ro{\advance\xdim by.5 \Einheit \advance\ydim by.5 \Einheit}%
  \def\rro{\advance\xdim by.25cm \advance\ydim by.5 \Einheit}%
  \def\roo{\advance\xdim by.5 \Einheit \advance\ydim by.25cm}%
  \def\l{\advance\xdim by-.30cm}%
  \def\r{\advance\xdim by.30cm}%
  \def\lu{\advance\xdim by-.5 \Einheit \advance\ydim by-.6 \Einheit}%
  \def\llu{\advance\xdim by-.25cm \advance\ydim by-.6 \Einheit}%
  \def\luu{\advance\xdim by-.5 \Einheit \advance\ydim by-.30cm}%
  \def\u{\advance\ydim by-.30cm}%
  \def\ru{\advance\xdim by.5 \Einheit \advance\ydim by-.6 \Einheit}%
  \def\rru{\advance\xdim by.25cm \advance\ydim by-.6 \Einheit}%
  \def\ruu{\advance\xdim by.5 \Einheit \advance\ydim by-.30cm}%
  #1\raise\ydim\hbox to0pt{\hskip\xdim
     \vbox to0pt{\vss\hbox to0pt{\hss$#2$\hss}\vss}\hss}%
}
\catcode`\@=13

\catcode`\@=11
\font@\twelverm=cmr10 scaled\magstep1
\font@\twelveit=cmti10 scaled\magstep1
\font@\twelvebf=cmbx10 scaled\magstep1
\font@\twelvei=cmmi10 scaled\magstep1
\font@\twelvesy=cmsy10 scaled\magstep1
\font@\twelveex=cmex10 scaled\magstep1

\newtoks\twelvepoint@
\def\twelvepoint{\normalbaselineskip15\p@
 \abovedisplayskip15\p@ plus3.6\p@ minus10.8\p@
 \belowdisplayskip\abovedisplayskip
 \abovedisplayshortskip\z@ plus3.6\p@
 \belowdisplayshortskip8.4\p@ plus3.6\p@ minus4.8\p@
 \textonlyfont@\rm\twelverm \textonlyfont@\it\twelveit
 \textonlyfont@\sl\twelvesl \textonlyfont@\bf\twelvebf
 \textonlyfont@\smc\twelvesmc \textonlyfont@\tt\twelvett
%Erg\"anzung des fetten Small-Capitals-Fonts:
%
 \ifsyntax@ \def\big##1{{\hbox{$\left##1\right.$}}}%
  \let\Big\big \let\bigg\big \let\Bigg\big
 \else
  \textfont\z@=\twelverm  \scriptfont\z@=\tenrm  \scriptscriptfont\z@=\sevenrm
  \textfont\@ne=\twelvei  \scriptfont\@ne=\teni  \scriptscriptfont\@ne=\seveni
  \textfont\tw@=\twelvesy \scriptfont\tw@=\tensy \scriptscriptfont\tw@=\sevensy
  \textfont\thr@@=\twelveex \scriptfont\thr@@=\tenex
        \scriptscriptfont\thr@@=\tenex
  \textfont\itfam=\twelveit \scriptfont\itfam=\tenit
        \scriptscriptfont\itfam=\tenit
  \textfont\bffam=\twelvebf \scriptfont\bffam=\tenbf
        \scriptscriptfont\bffam=\sevenbf
  \setbox\strutbox\hbox{\vrule height10.2\p@ depth4.2\p@ width\z@}%
  \setbox\strutbox@\hbox{\lower.6\normallineskiplimit\vbox{%
        \kern-\normallineskiplimit\copy\strutbox}}%
 \setbox\z@\vbox{\hbox{$($}\kern\z@}\bigsize@=1.4\ht\z@
 \fi
 \normalbaselines\rm\ex@.2326ex\jot3.6\ex@\the\twelvepoint@}

\font@\fourteenrm=cmr10 scaled\magstep2
\font@\fourteenit=cmti10 scaled\magstep2
\font@\fourteensl=cmsl10 scaled\magstep2
\font@\fourteensmc=cmcsc10 scaled\magstep2
\font@\fourteentt=cmtt10 scaled\magstep2
\font@\fourteenbf=cmbx10 scaled\magstep2
\font@\fourteeni=cmmi10 scaled\magstep2
\font@\fourteensy=cmsy10 scaled\magstep2
\font@\fourteenex=cmex10 scaled\magstep2
\font@\fourteenmsa=msam10 scaled\magstep2
\font@\fourteeneufm=eufm10 scaled\magstep2
\font@\fourteenmsb=msbm10 scaled\magstep2
\newtoks\fourteenpoint@
\def\fourteenpoint{\normalbaselineskip15\p@
 \abovedisplayskip18\p@ plus4.3\p@ minus12.9\p@
 \belowdisplayskip\abovedisplayskip
 \abovedisplayshortskip\z@ plus4.3\p@
 \belowdisplayshortskip10.1\p@ plus4.3\p@ minus5.8\p@
 \textonlyfont@\rm\fourteenrm \textonlyfont@\it\fourteenit
 \textonlyfont@\sl\fourteensl \textonlyfont@\bf\fourteenbf
 \textonlyfont@\smc\fourteensmc \textonlyfont@\tt\fourteentt
%Erg\"anzung des fetten Small-Capitals-Fonts:
%
 \ifsyntax@ \def\big##1{{\hbox{$\left##1\right.$}}}%
  \let\Big\big \let\bigg\big \let\Bigg\big
 \else
  \textfont\z@=\fourteenrm  \scriptfont\z@=\twelverm  \scriptscriptfont\z@=\tenrm
  \textfont\@ne=\fourteeni  \scriptfont\@ne=\twelvei  \scriptscriptfont\@ne=\teni
  \textfont\tw@=\fourteensy \scriptfont\tw@=\twelvesy \scriptscriptfont\tw@=\tensy
  \textfont\thr@@=\fourteenex \scriptfont\thr@@=\twelveex
        \scriptscriptfont\thr@@=\twelveex
  \textfont\itfam=\fourteenit \scriptfont\itfam=\twelveit
        \scriptscriptfont\itfam=\twelveit
  \textfont\bffam=\fourteenbf \scriptfont\bffam=\twelvebf
        \scriptscriptfont\bffam=\tenbf
  \setbox\strutbox\hbox{\vrule height12.2\p@ depth5\p@ width\z@}%
  \setbox\strutbox@\hbox{\lower.72\normallineskiplimit\vbox{%
        \kern-\normallineskiplimit\copy\strutbox}}%
 \setbox\z@\vbox{\hbox{$($}\kern\z@}\bigsize@=1.7\ht\z@
 \fi
 \normalbaselines\rm\ex@.2326ex\jot4.3\ex@\the\fourteenpoint@}

\font@\seventeenrm=cmr10 scaled\magstep3
\font@\seventeenit=cmti10 scaled\magstep3
\font@\seventeensl=cmsl10 scaled\magstep3
\font@\seventeensmc=cmcsc10 scaled\magstep3
\font@\seventeentt=cmtt10 scaled\magstep3
\font@\seventeenbf=cmbx10 scaled\magstep3
\font@\seventeeni=cmmi10 scaled\magstep3
\font@\seventeensy=cmsy10 scaled\magstep3
\font@\seventeenex=cmex10 scaled\magstep3
\font@\seventeenmsa=msam10 scaled\magstep3
\font@\seventeeneufm=eufm10 scaled\magstep3
\font@\seventeenmsb=msbm10 scaled\magstep3
\newtoks\seventeenpoint@
\def\seventeenpoint{\normalbaselineskip18\p@
 \abovedisplayskip21.6\p@ plus5.2\p@ minus15.4\p@
 \belowdisplayskip\abovedisplayskip
 \abovedisplayshortskip\z@ plus5.2\p@
 \belowdisplayshortskip12.1\p@ plus5.2\p@ minus7\p@
 \textonlyfont@\rm\seventeenrm \textonlyfont@\it\seventeenit
 \textonlyfont@\sl\seventeensl \textonlyfont@\bf\seventeenbf
 \textonlyfont@\smc\seventeensmc \textonlyfont@\tt\seventeentt
%Erg\"anzung des fetten Small-Capitals-Fonts:
%
 \ifsyntax@ \def\big##1{{\hbox{$\left##1\right.$}}}%
  \let\Big\big \let\bigg\big \let\Bigg\big
 \else
  \textfont\z@=\seventeenrm  \scriptfont\z@=\fourteenrm  \scriptscriptfont\z@=\twelverm
  \textfont\@ne=\seventeeni  \scriptfont\@ne=\fourteeni  \scriptscriptfont\@ne=\twelvei
  \textfont\tw@=\seventeensy \scriptfont\tw@=\fourteensy \scriptscriptfont\tw@=\twelvesy
  \textfont\thr@@=\seventeenex \scriptfont\thr@@=\fourteenex
        \scriptscriptfont\thr@@=\fourteenex
  \textfont\itfam=\seventeenit \scriptfont\itfam=\fourteenit
        \scriptscriptfont\itfam=\fourteenit
  \textfont\bffam=\seventeenbf \scriptfont\bffam=\fourteenbf
        \scriptscriptfont\bffam=\twelvebf
  \setbox\strutbox\hbox{\vrule height14.6\p@ depth6\p@ width\z@}%
  \setbox\strutbox@\hbox{\lower.86\normallineskiplimit\vbox{%
        \kern-\normallineskiplimit\copy\strutbox}}%
 \setbox\z@\vbox{\hbox{$($}\kern\z@}\bigsize@=2\ht\z@
 \fi
 \normalbaselines\rm\ex@.2326ex\jot5.2\ex@\the\seventeenpoint@}

\catcode`\@=13

\def\si{\sigma}
\def\la{\lambda}
\def\v#1{\vert#1\vert}
\def\({\left(}
\def\){\right)}

\def\VienAE{13}
\def\StanBI{12}
\def\StanAP{11}
\def\SiScAA{10}
\def\RoViAA{9}
\def\RoWZAA{8}
\def\RobeAA{7}
\def\OdlyAA{6}
\def\MaVaAA{5}
\def\KaPaAD{4}
\def\JaRiAA{3}
\def\FlajAA{2}
\def\ChWeAA{1}

\topmatter 
\title Permutations with restricted patterns and Dyck paths
\endtitle 
\author C.~Krattenthaler$^\dagger$
\endauthor 
\affil 
Institut f\"ur Mathematik der Universit\"at Wien,\\
Strudlhofgasse 4, A-1090 Wien, Austria.\\
e-mail: KRATT\@Ap.Univie.Ac.At\\
WWW: \tt http://radon.mat.univie.ac.at/People/kratt
\endaffil
\address Institut f\"ur Mathematik der Universit\"at Wien,
Strudlhofgasse 4, A-1090 Wien, Austria.
\endaddress
%\email KRATT\@Ap.Univie.Ac.At\\
%WWW: \tt http://radon.mat.univie.ac.at/People/kratt\endemail
%\dedicatory \enddedicatory
%\date \enddate
\thanks{$^\dagger$ Research partially supported by the Austrian
Science Foundation FWF, grant P13190-MAT}\endthanks
\subjclass Primary 05A05;
 Secondary 05A15, 05A16.
\endsubjclass
\keywords permutations with restricted patterns,
Chebyshev polynomials, continued fraction, Dyck paths\endkeywords
\abstract 
We exhibit a bijection between 132-avoiding permutations and Dyck
paths. Using this bijection, it is shown 
that all the recently discovered results on generating functions for
132-avoiding permutations with a given number of occurrences of the
pattern $12\dots k$ follow directly from old results on the
enumeration of Motzkin paths, among which is a continued fraction
result due to Flajolet.
As a bonus, we use these observations to derive further results and 
a precise asymptotic estimate for the number of 132-avoiding permutations of
$\{1,2,\dots,n\}$ with exactly $r$ occurrences of the pattern $12\dots k$.
Second, we exhibit a bijection between 123-avoiding permutations and
Dyck paths. When combined with a result of Roblet and Viennot, this
bijection allows us to express the generating function for
123-avoiding permutations with a given number of occurrences of 
the pattern $(k-1)(k-2)\dots 1k$ in form of a continued fraction and to
derive further results for these permutations.
\endabstract
\endtopmatter

\leftheadtext{C. Krattenthaler}

\document

\subhead 1. Introduction\endsubhead
In the recent papers \cite{\ChWeAA, \JaRiAA, \MaVaAA, \RoWZAA}, the authors
considered $132$-avoiding permutations with a prescribed number of
occurrences of the pattern $12\dots k$ (the most general results being
contained in \cite{\MaVaAA}) and $123$-avoiding permutations which
also avoid the pattern $(k-1)(k-2)\dots 1k$. 
(See the end of this section for the precise definition of
permutations which avoid a certain pattern, and of Dyck paths.)
They found that generating functions for
these permutations can be expressed in terms of continued fractions
and Chebyshev polynomials.

The purpose of this paper is to make a case for the paradigm:
\bigskip

{\leftskip1cm \rightskip1cm
%\twelvepoint
\noindent
\it `Whenever you encounter generating functions which can be
expressed in terms of continued fractions or Chebyshev polynomials,
then expect that Dyck or Motzkin paths are at the heart of
your problem, and will help to solve it.'
\par
}
\bigskip
Indeed, as I am going to demonstrate in Section~2, 
there is an obvious bijection
between $132$-avoiding permutations and Dyck paths. Known results for
generating functions for Motzkin paths (one of which
due to Flajolet \cite{\FlajAA}, the other being folklore; 
Dyck paths being special Motzkin
paths) then allow one immediately to express the
generating functions that we are interested in in terms of continued
fractions and Chebyshev polynomials (thus making the speculation in
\cite{\ChWeAA, Sec.~5} precise and explicit). In particular, we recover all
the relevant results from \cite{\ChWeAA, \MaVaAA, \RoWZAA}. Furthermore, by
exploiting the relation between $132$-avoiding permutations and Dyck
paths further, we are able to find an explicit expression for
the generating function for 132-avoiding permutations with exactly $r$
occurrences of the pattern $12\dots k$, thus extending a result from
\cite{\MaVaAA}. This, in turn, allows us to provide
a precise asymptotic estimate for the number of these permutations of
$\{1,2,\dots,n\}$ as $n$ becomes large. All these results can be found
in Section~3, as well as generating functions for
$132$-avoiding permutations with no occurrence or one occurrence of the
pattern $23\dots k1$.

In Section~4 we exhibit a bijection between $123$-avoiding
permutations and Dyck paths. In Section~5 we combine this bijection with a 
result of Roblet and Viennot \cite{\RoViAA} on the enumeration of Dyck
paths to obtain a continued fraction for the generating function of
$123$-avoiding permutations with a given number of occurrences of 
the pattern $(k-1)(k-2)\dots 1k$. Further results on these
permutations (which extend another result from \cite{\ChWeAA}) can be
found in Section~5 as well, including precise 
asymptotic estimates. (By combining this bijection between
$123$-avoiding permutations and Dyck paths with our bijection between
the latter and $132$-avoiding permutations, we obtain a bijection
between $123$-avoiding and $132$-avoiding permutations. This bijection
appears to be new. In particular, it is different from the one by
Simion and Schmidt \cite{\SiScAA, Sec.~6}, as can be immediately seen by
considering e.g\. the Examples on p.~404 of \cite{\SiScAA}.)

For the convenience of the reader, we recall the results on Motzkin
and Dyck paths, on which we rely so heavily, 
in an appendix at the end of the paper.

\medskip
At the end of the introduction, let us recall the basic definitions.

Let $\pi=\pi_1\pi_2\dots \pi_n$ be a permutation of $\{1,2,\dots,n\}$
and $\si=\si_1\si_2\dots\si_k$ be a permutation of $\{1,2,\dots,k\}$,
$k\le n$. We say that the permutation $\pi$ {\it contains the pattern}
$\si$, if there 
are indices $1\le i_1<i_2<\dots<i_k\le n$ such that 
$\pi_{i_1}\pi_{i_2}\dots \pi_{i_k}$ is in the same relative order as
$\si_1\si_2\dots\si_k$. Otherwise, $\pi$ is said to {\it avoid the pattern}
$\si$, or, alternatively, we say that $\pi$ is {\it $\si$-avoiding}.

A {\it Dyck path} is a lattice path in the plane integer lattice $\Bbb
Z^2$ ($\Bbb Z$ denoting the set of integers) consisting of up-steps
$(1,1)$ and down-steps $(1,-1)$, which never passes below the
$x$-axis. See Figure~1 for an example.

\subhead 2. A bijection between $132$-avoiding permutations and Dyck
paths\endsubhead
In this section we define a map $\Phi$ which maps $132$-avoiding
permutations to Dyck paths which start at the origin and return to the
$x$-axis as follows.
Let $\pi=\pi_1\pi_2\dots \pi_n$ be a $132$-avoiding permutation. 
We read the permutation $\pi$ from left to right and successively
generate a Dyck path. When $\pi_j$ is read, then in the path we adjoin
as many up-steps as necessary, followed by a down-step from height
$h_j+1$ to height $h_j$ (measured from the $x$-axis), where $h_j$ is the
number of elements in $\pi_{j+1}\dots\pi_n$ which are larger than
$\pi_j$. 

\midinsert
\vskip10pt
\vbox{
$$
\Gitter(17,6)(0,0)
\Koordinatenachsen(17,6)(0,0)
\Pfad(0,0),3343334344344434\endPfad
\DickPunkt(0,0)
\DickPunkt(16,0)
\hskip9cm
$$
\centerline{\eightpoint The Dyck path corresponding to $74352681$}
\vskip7pt
\centerline{\eightpoint Figure 1}
}
\vskip10pt
\endinsert

For example, let $\pi=74352681$. The first element to be read is
$7$. There is $1$ element in $4352681$ which is larger than $7$,
therefore the path starts with two up-steps followed by a down-step,
thus reaching height $1$ (see Figure~1). Next $4$ is read. There are
$3$ elements in $352681$ which are larger than $4$, therefore the
path continues with three up-steps followed by a down-step, thus
reaching height $3$. Etc. The complete Dyck path
$\Phi(74352681)$ is shown in Figure~1.

The reader should note that, for the map $\Phi$ to be well-defined, it
is essential that the permutation $\pi$ to which the map is applied is
$132$-avoiding. For this guarantees that always $h_j-1\le h_{j+1}$, so
that it is always possible to connect the down-step from height $h_j+1$
to $h_j$ (formed by definition of $\Phi$ when considering $\pi_j$) by
a number of up-steps (this number being possibly zero) to the 
down-step from height $h_{j+1}+1$ to $h_{j+1}$. Conversely, given a
Dyck path starting at the origin and returning to the $x$-axis, the
obvious inverse of $\Phi$ produces a $132$-avoiding permutation.

In summary, the map $\Phi$ is a bijection between $132$-avoiding
permutations of\break $\{1,2,\dots,n\}$ and Dyck paths from $(0,0)$ to
$(2n,0)$. We remark that, in view of the standard bijection between
rooted ordered trees and Dyck paths through a depth-first traversal of
the trees (cf\. e.g\. \cite{\StanBI, Prop.~6.2.1 (i) and (v),
Cor.~6.2.3 (i) and (v)}), this map is equivalent to a bijection between 
$132$-avoiding permutations and rooted ordered trees given by Jani and
Rieper \cite{\JaRiAA}.

For the sake of completeness, and to show the close relation between 
the map $\Phi$ and the map $\Psi$ that is to be defined in Section~4,
we provide an alternative way to define the map $\Phi$. 
Let $\pi=\pi_1\pi_2\dots \pi_n$ be a $132$-avoiding permutation. 
In $\pi$, we determine all the {\it left-to-right minima}. A left-to-right
minimum is an element $\pi_i$ which is smaller than all the elements to
its left, i.e., smaller than all $\pi_j$ with $j<i$. For example, 
the left-to-right minima in the permutation $74352681$ are
$7$, $4$, $3$, $2$, $1$. 

Let the left-to-right minima in $\pi$ be $m_1$, $m_2$, \dots, $m_s$,
so that
$$\pi=m_1w_1m_2w_2\dots m_sw_s,\tag3.1$$
where $w_i$ is the subword of $\pi$ in between $m_{i}$ and
$m_{i+1}$. Read the decomposition (3.1)
from left to right. Any left-to-right minimum $m_i$ is translated into
$m_{i-1}-m_{i}$ up-steps (with the convention $m_{0}=n+1$). Any subword
$w_i$ is translated into $\v{w_i}+1$ down-steps (where
$\v{w_i}$ denotes the number of elements of $w_i$). 

\medskip
In the lemma below we list two properties of the bijection $\Phi$, which
will be subsequently used in Section~3.

\proclaim{Lemma $\Phi$}
Let $\pi=\pi_1\pi_2\dots\pi_n$ be a $132$-avoiding permutation, and let
$P=\Phi(\pi)$ be the corresponding Dyck path. Then, 
\roster 
\item a down-step in $P$ from height $i$ to height $i-1$ corresponds in
a one-to-one fashion to an element $\pi_j$ in the permutation which is
the first element in an increasing subsequence in $\pi$ of
length $i$ (i.e., an occurrence of the pattern $12\dots i$) 
that is maximal with respect to the property that $\pi_j$
is its first element.
\item a portion of the path $P$ which starts at height $h+i-1$ and
eventually falls down to height $h$, for some $h$, 
followed by an up-step corresponds
to an occurrence of the pattern $23\dots i1$ in $\pi$.
\endroster
\endproclaim

\demo{Proof}
Re (1): By the definition of $\Phi$, a down-step from height $i$ to
height $i-1$ in $P=\Phi(\pi)$ means that we read an element $\pi_j$
which has the property that there are $i-1$ elements in
$\pi_{j+1}\dots \pi_n$ which are larger than $\pi_j$. Since $\pi$ is
$132$-avoiding, these $i-1$ elements have to appear in increasing order,
thus, together with $ \pi_j$, form an increasing subsequence of
length $i$ that cannot be made longer under the assumption that $\pi_j$
is the first element in the subsequence.

Re (2): In a path portion which starts at height $h+i-1$ and eventually
falls down to height $i$ we must find a down-step from height $h+i-1$ to
$h+i-2$, a down-step from height $h+i-2$ to
$h+i-3$, \dots, a down-step from height $i+1$ to $i$. Under the
correspondence $\Phi$, these down-steps correspond to an increasing 
subsequence $\pi_{j_1}\pi_{j_2}\dots\pi_{j_{i-1}}$ of length $i-1$ in
$\pi$. If now the path continues by (at least) one up-step, 
then the following down-step corresponds to an element $\pi_{j_i}$,
$j_i>j_{i-1}$,
with the property that there are more elements in
$\pi_{j_{i+1}}\dots\pi_n$ that are larger than $\pi_{j_i}$ than there
are elements in $\pi_{j_{i}}\dots\pi_n$ that are larger than
$\pi_{j_{i-1}}$. Evidently, this is only possible if
$\pi_{j_{i-1}}>\pi_{j_{i}}$. Since $\pi$ is $132$-avoiding, this implies
that we have even $\pi_{j_1}>\pi_{j_{i}}$. Hence, 
$\pi_{j_1}\pi_{j_2}\dots\pi_{j_{i-1}}\pi_{j_i}$ is an occurrence of the
pattern $23\dots i1$ in $\pi$. \qed
\enddemo

\subhead 3. The enumeration of $132$-avoiding permutations with a
prescribed number of occurrences of the patterns $12\dots k$ and
$23\dots k1$\endsubhead
In this section we provide explicit expressions for generating
functions for $132$-avoiding permutations with a prescribed number of
occurrences of the pattern $12\dots k$, and for $132$-avoiding
permutations with a prescribed number of occurrences of the pattern
$23\dots k1$. 

First we consider the former permutations. Given a $132$-avoiding
permutation $\pi$, we denote the number of occurrences of the pattern
$12\dots k$ in $\pi$ by $N(12\dots k;\pi)$. 

Given a Dyck path $P$, we assign a weight to it, denoted by $w_1(k;P)$.
It is defined as the sum $\sum _{d} ^{}\binom {i(d)-1}{k-1}$, where
the sum is over all down-steps $d$ of $P$, and where $i(d)$ is the
height of the starting point of $d$. For example, the weight
$w_1(k;.)$ of the Dyck path in Figure~1 is
$$
\binom12+\binom32+\binom32+\binom22+\binom22+\binom12+\binom02+\binom02=8.
$$
From Lemma~$\Phi$.(1) it is immediate that
$$N(12\dots k;\pi)=w_1(k;\Phi(\pi)).\tag3.2$$
This observation, combined with Flajolet's continued fraction theorem
for the generating function of Motzkin paths (see Theorem~A1), allows us
to express the generating function which counts $132$-avoiding permutations 
with respect to the number of occurrences of the pattern $12\dots k$ in
form of a continued fraction. This result was first obtained by
Mansour and Vainshtein \cite{\MaVaAA, Theorem~2.1}.
In the statement of the theorem, and in the following, we write
$\v{\pi}$ for the number of elements which are permuted by $\pi$. For
example, we have $\v{74352681}=8$.

\proclaim{Theorem 1}
The generating function $\sum _{\pi} ^{}y^{N(12\dots
k;\pi)}x^{\v{\pi}}$, where
the sum is over all $132$-avoiding permutations, is given by 
$$
\cfrac 1\\1-
\cfrac xy^{\binom 0{k-1}}\\1-
\cfrac xy^{\binom 1{k-1}}\\1-
\cfrac xy^{\binom 2{k-1}}\\1-
\cdots\endcfrac\endcfrac\endcfrac\endcfrac\ .
\tag 3.3$$
\endproclaim
\demo{Proof} Apply Theorem~A1 with $b_i=0$ and
$\la_i=xq^{\binom{i-1}{k-1}}$, $i=0,1,\dots$, and use (3.2). \qed
\enddemo

We remark that the above proof is essentially equivalent to the one in
\cite{\JaRiAA, proof of Corollary~7}. It is obvious that the
refinement in \cite{\MaVaAA, expression for $W_k(\dots)$ after
Proposition~2.3} could also easily be derived by using the
correspondence $\Phi$ and Flajolet's continued fraction.

\medskip
Next we turn our attention to $132$-avoiding permutations with a {\it
fixed\/} 
number of occurrences of the pattern $12\dots k$. The theorem below was
first obtained by Chow and West \cite{\ChWeAA, Theorem~3.6, second
case} in an equivalent form.

\proclaim{Theorem 2}
The generating function $\sum _{\pi} ^{}x^{\v{\pi}}$, where
the sum is over all $132$-avoiding permutations 
which also avoid the pattern $12\dots k$, is given by
$$\frac {U_{k-1}\(\frac {1} {2\sqrt x}\)} {\sqrt xU_{k}\(\frac {1}
{2\sqrt x}\)},\tag3.4$$
where $U_n(x)$ denotes the $n$-th Chebyshev polynomial of the second
kind, $U_n(\cos t)=\sin((n+1)t)/\sin t$.
\endproclaim
\demo{Proof} By Lemma~$\Phi$.(1), the permutations in the statement of
the theorem are in bijection with Dyck paths, which start at the
origin, return to the $x$-axis, and do not exceed the height $k-1$. Now apply
Theorem~A2 with $b_i=0$, $\la_i=1$, $i=0,1,\dots$, $K=k-1$, $r=s=0$, $x$
replaced by $\sqrt x$, and use Fact~A3. \qed
\enddemo

The next theorem extends a result by Mansour and Vainshtein
\cite{\MaVaAA, Theorems~3.1 and 4.1}, who proved the special case when
$r$ is at most $k(k+3)/2$.

\proclaim{Theorem 3}
Let $r\ge1$.
The generating function $\sum _{\pi} ^{}x^{\v{\pi}}$, where
the sum is over all $132$-avoiding permutations 
with exactly $r$
occurrences of the pattern $12\dots k$, is given by
$$\sum _{} ^{}\(\binom {\ell_1+\ell_2-1}{\ell_2}
\binom {\ell_2+\ell_3-1}{\ell_3}\cdots\)
\frac {\(U_{k-1}\(\frac {1} {2\sqrt x}\)\)^{\ell_1-1}} 
{\(U_{k}\(\frac {1} {2\sqrt x}\)\)^{\ell_1+1}}
x^{\frac {1} {2}(\ell_1-1)+(\ell_2+\ell_3+\cdots)},
\tag3.5$$
where the sum is over all nonnegative integers $\ell_1,\ell_2,\dots$
with
$$\ell_1\binom{k-1}{k-1}+\ell_2\binom{k}{k-1}+\ell_3\binom{k+1}{k-1}+\cdots
=r,\tag3.6$$
and where $U_n(x)$ denotes the $n$-th Chebyshev polynomial of the second
kind.
\endproclaim

\remark{Remark} Because of (3.6), almost all summation indices of the
sum in (3.5) must be zero, so that the sum in (3.5) is in fact a
finite sum. In particular, it reduces to just one term if $r\le k$,
thus recovering \cite{\MaVaAA, Theorem~3.1}, and it reduces to a
single sum if $k<r\le k(k+3)/2$, thus recovering \cite{\MaVaAA,
Theorem~4.1}. 

\endremark

\demo{Proof of Theorem 3}
Let $\pi$ be a permutation of the statement of the theorem. We apply
$\Phi$ to obtain the corresponding Dyck path $P=\Phi(\pi)$. The Dyck
path $P$ has a unique decomposition of the form
$$P_0V_1d_1P_1V_2d_2P_2\dots V_sd_sP_s,\tag3.7$$
where $P_0$ is the portion of $P$ from the origin until the first time
the height $k-1$ is reached, where the $d_i$'s are the down steps whose end
points have at least the height $k-1$, where the $P_i$'s,
$i=1,2,\dots,s-1$, are path portions which start and end at height
$k-1$ and never exceed height $k-1$, where the $V_i$'s are
path portions consisting of several subsequent up-steps which fill
the gaps in between, and where 
$P_s$ is the portion of $P$ from the last point at height $k-1$ until
the end of the path. The path portion $P_i$ can only be nonempty if
$d_i$ is a down-step from height $k$ to
$k-1$. Clearly, $d_s$ {\it must be} a down-step from height $k$ to
height $k-1$.  

Now suppose that among the $d_i$'s there are $\ell_1$ down-steps from
height $k$ to height $k-1$, $\ell_2$ down-steps from
height $k+1$ to height $k$, etc. Because of (3.2) the relation (3.6)
must hold. 

Let us for the moment fix $\ell_1,\ell_2,\dots$ and ask how many 
orderings of $\ell_1$ down-steps from
height $k$ to height $k-1$, $\ell_2$ down-steps from
height $k+1$ to height $k$, etc., there are
which can come from a decomposition
(3.7) when we ignore the $P_i$'s and $V_i$'s. In fact, there are many
restrictions to be obeyed: After a 
down-step from height $h+1$ to height $h$ there can only follow a down
step of at least that height or at worst from height $h$ to height
$h-1$. Let $t$ be maximal so that $\ell_t$ is nonzero. Then the above
observation 
tells that after a down-step from height $t+k-1$ to $t+k-2$ there can
only follow another down-step of this sort (of course, with an up-step
in between) or a down-step from height $t+k-2$ to $t+k-3$. Hence, if
we just concentrate on these two types of down-steps, of which there
are $\ell_t$ and $\ell_{t-1}$, respectively, then there are exactly $\binom
{\ell_{t-1}+\ell_t-1}{\ell_t}$ different possible orderings between
these steps, taking into account that the last step out of these must
necessarily be a step from height $t+k-2$ to $t+k-3$. Next, by similar
considerations, one concludes that, given an ordering of the
down-steps from height $t+k-1$ to $t+k-2$ and from height $t+k-2$ to
$t+k-3$, there are $\binom
{\ell_{t-2}+\ell_{t-1}-1}{\ell_{t-1}}$ different possibilities to
intersperse $\ell_{t-2}$ down-steps from height $t+k-3$ to
$t+k-4$. Etc. This explains the product of binomials in (3.5). 

To explain the remaining expression, we observe that after any of the 
$\ell_1$ down-steps, from height $k$ to height $k-1$, $d_i$ say,
except for the last, there follows a (possibly empty) path $P_i$,
which is a Dyck path which starts and ends at height $k-1$ and never
exceeds height $k-1$. By Theorem~A2 with $b_i=0$, $\la_i=1$,
$i=0,1,\dots$, $K=k-1$, $r=s=k-1$, $x$ 
replaced by $\sqrt x$, and use of Fact~A3, we conclude that 
the generating function for these paths is equal to
$$\frac {x^{(k-1)/2}U_{k-1}\(\frac {1} {2\sqrt x}\)} 
{x^{k/2}U_{k}\(\frac {1} {2\sqrt x}\)}.$$
By again applying Theorem~A2, this time with $b_i=0$, $\la_i=1$,
$i=0,1,\dots$, $K=k-1$, $r=0$, $s=k-1$, $x$ 
replaced by $\sqrt x$, and using Fact~A3, we obtain that
the generating function for paths $P_0$ from the origin to
height $k-1$, never exceeding height $k-1$, is equal to
$$\frac {x^{(k-1)/2}} {x^{k/2}U_{k}\(\frac {1} {2\sqrt x}\)}.\tag3.8$$
Similarly, the generating function for paths $P_s$ from height $k-1$ to
height $0$, never exceeding height $k-1$, is also given by (3.8). If
everything is combined, the expression (3.5) results.
\qed
\enddemo

Theorem~3 can be readily used to find an asymptotic formula for the
number of $132$-avoiding permutations with exactly $r$ occurrences of
the pattern $12\dots k$. The corresponding result, given in the
theorem below, extends \cite{\ChWeAA, Corollary~4.2}.
Before we state the theorem, we recall an elementary lemma
(cf\. e.g\. \cite{\OdlyAA, Sec.~9.1}).

\proclaim{Lemma 4}
Let $f(x)$ and $g(x)$ be polynomials. It is assumed that all the
zeroes of $g(x)$ have modulus larger than $a>0$.
Consider the expansion
$$\frac {f(x)} {(x-a)^Rg(x)}=\sum _{n=0} ^{\infty}c_nx^n.$$
Then, as $n$ becomes large, we have
$$c_n\sim (-1)^R\frac {n^{R-1}} {(R-1)!}a^{-n-R}\frac {f(a)}
{g(a)}\(1+O\(\frac {1} {n}\)\).\eqno \qed$$
\endproclaim

\proclaim{Theorem 5}
Let $r$ and $k$ be fixed nonnegative integers. Then, as $n$ becomes
large, the number of $132$-avoiding permutations with exactly $r$
occurrences of the pattern $12\dots k$ is asymptotically
$$\(\frac {4\sin^2\frac {\pi} {k+1}} {k+1}\)^{r+1}\frac {n^r} {r!}
\(4\cos^2\frac {\pi} {k+1}\)^{n-r}\(1+O\(\frac {1} {n}\)\).\tag3.9$$ 
\endproclaim
\demo{Proof}If $r=0$, this follows immediately from Theorem~2 
and Lemma~4. If $r\ge1$, we start from the generating
function given in Theorem~3. It should be observed that, in view of
Lemma~4, the summand in (3.5) which asymptotically yields the largest
contribution is the one with $\ell_1=r$ and all other $\ell_i$'s equal
to zero. Then application of Lemma~4 to this summand gives (3.9) after
some computation.
\qed
\enddemo

The next group of results concerns the enumeration of $132$-avoiding
permutations with a given number of occurrences of the pattern $23\dots
k1$. We use again the map $\Phi$ to translate these permutations 
into Dyck paths. The property of $\Phi$ which is important now is
given by Lemma~$\Phi$.(2). It says that we can recognize the occurrence
of a pattern $23\dots k1$ in a $132$-avoiding permutation in the
corresponding Dyck path by a portion of the path 
which starts at height $h+k-1$, 
eventually falls down to height $k$, and is then followed by an
up-step.

As the first application of our approach we show how to rederive
another result due to Chow and West \cite{\ChWeAA, Theorem~3.6, third
case}, which is reformulated here in an equivalent form.

\proclaim{Theorem 6}
The generating function $\sum _{\pi} ^{}x^{\v{\pi}}$, where
the sum is over all $132$-avoiding permutations 
which also avoid the pattern $23\dots k1$, is given by
$$\frac {U_{k-1}\(\frac {1} {2\sqrt x}\)} {\sqrt xU_{k}\(\frac {1}
{2\sqrt x}\)},\tag3.10$$
where $U_n(x)$ denotes the $n$-th Chebyshev polynomial of the second
kind.
\endproclaim
\demo{Proof}
Let $\pi$ be a permutation of the statement of the theorem.
By the observation above the statement of the theorem (which was 
based on Lemma~$\Phi$.(2)), for the corresponding Dyck path
$P=\Phi(\pi)$ there are two possibilities: Either $P$ never exceeds
the height $k-2$ (and, thus, $\pi$ does not contain any increasing
subsequence of length $k-1$), or $P$ can be decomposed as
$$P_0u_0P_1u_1\dots P_{s}u_{s}P_{s+1}D,\tag3.11$$
where $P_0$ is a path from the origin to height $k-2$ never exceeding
height $k-2$,
where for $i=1,2,\dots,s$ the portion $P_i$ is a path starting and
ending at height $i+k-2$, never running below height $i$,
and never exceeding height $i+k-2$, where for $i=0,1,\dots,s$ the step
$u_i$ is an up-step from height $i+k-2$ to height $i+k-1$, where
$P_{s+1}$ is a path from height $s+k-1$ to $s+1$, never running below
height $s+1$, and never exceeding height $s+k-1$, and where $D$
consists of $s+1$ down-steps, from height $s+1$ to height $0$. 

By Theorem~A2 with $b_i=0$, $\la_i=1$,
$i=0,1,\dots$, $K=k-2$, $r=s=0$, $x$ 
replaced by $\sqrt x$, and Fact~A3, the generating function for the
Dyck paths which never exceed height $k-2$ is equal to
$$\frac {x^{(k-2)/2}U_{k-2}\(\frac {1} {2\sqrt x}\)} 
{x^{(k-1)/2}U_{k-1}\(\frac {1} {2\sqrt x}\)}.\tag3.12$$
By Theorem~A2 with $b_i=0$, $\la_i=1$,
$i=0,1,\dots$, $K=k-2$, $r=0$, $s=k-2$, $x$ 
replaced by $\sqrt x$, and Fact~A3, the generating function for the
possible paths $P_0$ in the decomposition (3.11) is equal to
$$\frac {x^{(k-2)/2}} 
{x^{(k-1)/2}U_{k-1}\(\frac {1} {2\sqrt x}\)},$$
as well as the generating function for the possible paths
$P_{s+1}$. Finally, for any fixed $j$ between $1$ and $s$, 
by Theorem~A2 with $b_i=0$, $\la_i=1$,
$i=0,1,\dots$, $K=k-2$, $r=s=k-2$, $x$ 
replaced by $\sqrt x$, and Fact~A3, the generating function for the
possible paths $P_j$ in the decomposition (3.11) is also given by
(3.12). 

If everything is combined, then we obtain that the generating function
for the permutations of the statement of the theorem is given by
$$\frac {U_{k-2}\(\frac {1} {2\sqrt x}\)} 
{\sqrt xU_{k-1}\(\frac {1} {2\sqrt x}\)}+
\sum _{s\ge0} ^{}\frac {1} 
{\sqrt xU_{k-1}\(\frac {1} {2\sqrt x}\)}
\(\frac {U_{k-2}\(\frac {1} {2\sqrt x}\)} 
{\sqrt xU_{k-1}\(\frac {1} {2\sqrt x}\)}\)^s
\frac {1} 
{\sqrt xU_{k-1}\(\frac {1} {2\sqrt x}\)}x^{s+1}.$$
The sum is a geometric series and can therefore be evaluated. It is
then routine to convert the resulting expression 
into the expression (3.10), by
using standard identities for the Chebyshev polynomials. \qed 
\enddemo

It seems difficult to find an explicit expression for the generating
function for $132$-avoiding permutations 
with exactly $r$ occurrences of the pattern $23\dots k1$ for general
$r$. Yet, as long as
$1\le r\le k-1$ such an explicit expression can be easily derived.

\proclaim{Theorem 7}
The generating function $\sum _{\pi} ^{}x^{\v{\pi}}$, where
the sum is over all $132$-avoiding permutations 
with exactly one occurrence of the pattern $23\dots k1$, is given by
$$\frac {x} {U_{k-2}\(\frac {1}
{2\sqrt x}\)U_{k}\(\frac {1}
{2\sqrt x}\)},\tag3.13$$
where $U_n(x)$ denotes the $n$-th Chebyshev polynomial of the second
kind.

More generally, let $1\le r\le k-1$. Then
the generating function $\sum _{\pi} ^{}x^{\v{\pi}}$, where
the sum is over all $132$-avoiding permutations 
with exactly $r$ occurrences of the pattern $23\dots k1$, is given by
$$\frac {1} {U_{k-3}\(\frac {1}
{2\sqrt x}\)U_{k}\(\frac {1}
{2\sqrt x}\)}\sum _{\ell\mid r} ^{}\frac {1} {\ell+1}\binom
{2\ell}\ell x^{\ell+\frac {r} {2\ell}-\frac {1} {2}}
\(\frac {U_{k-3}\(\frac {1} {2\sqrt x}\)} 
{U_{k-2}\(\frac {1} {2\sqrt x}\)}\)^{r/\ell}.\tag3.14$$
\endproclaim
\demo{Proof}
Let $\pi$ be a $132$-avoiding permutation
with exactly $r$ occurrences of the pattern $23\dots k1$.
Analogously to the argument in the proof of Theorem~6, the Dyck path
$\Phi(\pi)$ corresponding to $\pi$ can be decomposed as
$$P_0u_0P_1u_1\dots P_{s-1}u_{s-1}P_{s}u_sd_s
P_{s+1}u_{s+1}d_{s+1}P_{s+2}
\dots u_{s+\ell-1}d_{s+\ell-1}P_{s+\ell}d_{s+\ell}P_{s+\ell+1}D,$$
where for $i=0,1,\dots,s$ the portions $P_i$ and steps $u_i$ are as in
(3.11), where the step $d_{s}$
is a down-step from height $s+k-1$ to height
$s+k-2$, where for $i=1,2,\dots,\ell-1$ the path $P_{s+i}$ is a path
starting and ending at height $s+k-2$, never running below 
height $s+1$, and never exceeding height $s+k-2$, 
the step $u_{s+i}$ is an
up-step from height $s+k-2$ to height $s+k-1$, and the step $d_{s+i}$
is a down-step from height $s+k-1$ to height $s+k-2$, where
$P_{s+\ell}$ is a path from height $s+k-2$ to $s+1$, never running below
height $s+1$, and never exceeding height $s+k-2$, where $d_{s+\ell}$ is a
down-step from height $s+1$ to height $s$, where $P_{s+\ell+1}$ is a path
of length $2r/\ell$,
starting and ending at height $s$, and never running below height $s$,
and where $D$ 
consists of $s$ down-steps, from height $s$ to height $0$. 

Still following the arguments in the proof of Theorem~6, and taking
into account that the number of possible paths $P_{s+\ell+1}$ is the
$(r/\ell)$-th Catalan number, this
decomposition implies that the generating function that we are looking
for is given by
$$\multline
\sum _{\ell\mid r} ^{}\sum _{s\ge0} ^{}\frac {1} 
{\sqrt xU_{k-1}\(\frac {1} {2\sqrt x}\)}
\(\frac {U_{k-2}\(\frac {1} {2\sqrt x}\)} 
{\sqrt xU_{k-1}\(\frac {1} {2\sqrt x}\)}\)^s\\
\cdot\(\frac {U_{k-3}\(\frac {1} {2\sqrt x}\)} 
{\sqrt xU_{k-2}\(\frac {1} {2\sqrt x}\)}\)^{\ell-1}
\frac {1} 
{\sqrt xU_{k-2}\(\frac {1} {2\sqrt x}\)}
\frac {1} {r/\ell+1}\binom
{2r/\ell}{r/\ell} x^{s+\ell+\frac {r} {\ell}+1/2}.
\endmultline$$
The inner sum is again a geometric series and can therefore be
evaluated. A routine calculation, followed by a replacement of $\ell$
by $r/\ell$, then transforms the resulting expression into (3.13).
\qed
\enddemo

It is obvious that in both cases (that is, for $132$-avoiding
permutations with no occurrence of the pattern $23\dots k1$,
respectively with $r\le k-1$ occurrences) Lemma~4 could be applied to
derive asymptotic formulas for the number of such permutations of
$\{1,2,\dots,n\}$, as $n$ becomes large. We omit the statement of the
corresponding formulas for the sake of brevity.  

It appears that, for generic $r$ (i.e., also for $r\ge k$), 
the number of $132$-avoiding
permutations of $\{1,2,\dots, n\}$ with exactly $r$
occurrences of the 
pattern $23\dots k1$ is asymptotically of the order
$\Theta\Big((4\cos^{2} \frac {\pi} {k+1})^n\Big)$, but we are not able
to offer a rigorous proof.

\subhead 4. A bijection between $123$-avoiding permutations and Dyck
paths\endsubhead
In this section we define another map, $\Psi$, between permutations
and Dyck paths, which maps $123$-avoiding
permutations to Dyck paths which start in the origin and return to the
$x$-axis.

Let $\pi=\pi_1\pi_2\dots \pi_n$ be a $123$-avoiding permutation. 
In $\pi$, we determine all the {\it right-to-left maxima}. A right-to-left
maximum is an element $\pi_i$ which is larger than all the elements to
its right, i.e., larger than all $\pi_j$ with $j>i$. For example, 
the right-to-left maxima in the permutation $58327641$ are
$1$, $4$, $6$, $7$, $8$.

Let the right-to-left maxima in $\pi$ be $m_1$, $m_2$, \dots, $m_s$,
from right to left, so that
$$\pi=w_sm_sw_{s-1}m_{s-1}\dots w_1m_1,\tag4.1$$
where $w_i$ is the subword of $\pi$ in between $m_{i+1}$ and
$m_i$. Since $\pi$ is $123$-avoiding, for all $i$ the elements in
$w_i$ must be in decreasing order. Moreover, for all $i$ all the
elements of $w_i$ are smaller than all the elements of $w_{i+1}$.

Now we are able to define the map $\Psi$. Read the decomposition (4.1)
from right to left. Any right-to-left maximum $m_i$ is translated into
$m_i-m_{i-1}$ up-steps (with the convention $m_{0}=0$). Any subword
$w_i$ is translated into $\v{w_i}+1$ down-steps (where, again,
$\v{w_i}$ denotes the number of elements of $w_i$). Finally, the
resulting path is reflected into a vertical line. (Alternatively, we
could have said that we generate the Dyck path from the back to the front.)
The Dyck path which corresponds to our special permutation $58327641$
is the one in Figure~1.

It is easy to see that the map $\Psi$ is a bijection between $123$-avoiding
permutations of $\{1,2,\dots,n\}$ and Dyck paths from $(0,0)$ to
$(2n,0)$.
The lemma below states the crucial property of this bijection, which
will be subsequently used in Section~5.

\proclaim{Lemma $\Psi$}
Let $\pi=\pi_1\pi_2\dots\pi_n$ be a $123$-avoiding permutation and let
$P=\Psi(\pi)$ be the corresponding Dyck path. Then
a {\rm peak} in $P$ of height $i$ (i.e., an up-step from height $i-1$ to
height $i$ followed by a down-step from height $i$ to height $i-1$) 
corresponds in
a one-to-one fashion to an element $\pi_j$ in the permutation which is
the last element in an occurrence of the pattern $(i-1)(i-2)\dots 1i$
that is maximal with respect to the property that $\pi_j$
is its last element.
\endproclaim

\demo{Proof} By construction of $\Psi$, any peak in the Dyck path
corresponds to a right-to-left maximum, $m$ say, in the
permutation. Furthermore, by induction one sees that the height of the
peak is equal to the number of elements to the left of $m$ that are
smaller than $m$. Clearly, all these elements belong to some $w_j$ in
the decomposition (4.1) of the permutation. By the observations above,
these elements are in decreasing order, and, thus, together with $m$
form an occurrence of the pattern $(i-1)(i-2)\dots 1i$
that cannot be made longer under the assumption that $m$
is the last element in the occurrence of the pattern. 
This proves the assertion of
the lemma. \qed
\enddemo

\subhead 5. The enumeration of $123$-avoiding permutations with a
prescribed number of occurrences of the pattern
$(k-1)(k-2)\dots 1k$\endsubhead
Let $\pi$ be a $132$-avoiding
permutation. We denote the number of occurrences of the pattern
$(k-1)(k-2)\dots 1k$ in $\pi$ by $N((k-1)(k-2)\dots 1k;\pi)$.

Given a Dyck path $P$, we assign a weight to it, denoted by $w_2(k;P)$.
It is defined as the sum $\sum _{d} ^{}\binom {i(p)-1}{k-1}$, where
the sum is over all peaks $p$ of $P$,
and where $i(p)$ is the
height of the peak. For example, the weight
$w_2(k;.)$ of the Dyck path in Figure~1 is
$$
\binom12+\binom32+\binom32+\binom22+\binom02=7.
$$
From Lemma~$\Psi$ it is immediate that
$$N((k-1)(k-2)\dots 1k;\pi)=w_1(k;\Psi(\pi)).\tag4.2$$
This observation, combined with Roblet and Viennot's 
continued fraction theorem
for the generating function of Dyck paths (see Theorem~A5), allows us
to express the generating function which counts $123$-avoiding permutations 
with respect to the number of occurrences of the pattern
$(k-1)(k-2)\dots 1k$ in
form of a continued fraction. 
Again, in the statement of the theorem, we write
$\v{\pi}$ for the number of elements which are permuted by $\pi$. 

\proclaim{Theorem 8}
The generating function $\sum _{\pi} ^{}y^{N((k-1)(k-2)\dots 1k;\pi)}
x^{\v{\pi}}$, where
the sum is over all $123$-avoiding permutations, is given by 
$$
\cfrac 1\\1-x\(y^{\binom 0{k-1}}-1\)-
\cfrac x\\1-x\(y^{\binom 1{k-1}}-1\)-
\cfrac x\\1-x\(y^{\binom 2{k-1}}-1\)-
\cdots\endcfrac\endcfrac\endcfrac\ .
\tag 4.3$$
\endproclaim
\demo{Proof} Apply Theorem~A5 with 
$\la_i=x$ and $\nu_i=xy^{\binom{i-1}{k-1}}$, $i=0,1,\dots$, and use
(4.2). \qed
\enddemo

Next, similar to Section~3, we study generating functions for
$123$-avoiding permutations with a {\it fixed\/} number of occurrences
of the pattern $(k-1)(k-2)\dots 1k$. The first theorem restates a
result due to Chow and West \cite{\ChWeAA, Theorem~3.6, first
case}. The proof however is different, as it is based on 
our Dyck path approach.

\proclaim{Theorem 9}
The generating function $\sum _{\pi} ^{}x^{\v{\pi}}$, where
the sum is over all $123$-avoiding permutations 
which also avoid the pattern $(k-1)(k-2)\dots 1k$, is given by
$$\frac {U_{k-1}\(\frac {1} {2\sqrt x}\)} {\sqrt xU_{k}\(\frac {1}
{2\sqrt x}\)},\tag4.4$$
where $U_n(x)$ denotes the $n$-th Chebyshev polynomial of the second
kind.
\endproclaim
\demo{Proof} By Lemma~$\Psi$, the permutations in the statement of
the theorem are in bijection with Dyck paths, which start at the
origin, return to the $x$-axis, and do not exceed the height
$k-1$. Now we apply
Theorem~A2 with $b_i=0$, $\la_i=1$, $i=0,1,\dots$, $K=k-1$, $r=s=0$, $x$
replaced by $\sqrt x$, and use Fact~A3. \qed
\enddemo

\proclaim{Theorem 10}
Let $1\le r\le k-1$.
The generating function $\sum _{\pi} ^{}x^{\v{\pi}}$, where
the sum is over all $123$-avoiding permutations 
with exactly $r$
occurrences of the pattern $(k-1)(k-2)\dots 1k$, is given by
$$x^{(r-1)/2}\frac {\(U_{k-1}\(\frac {1} {2\sqrt x}\)\)^{r-1}} 
{\(U_{k}\(\frac {1} {2\sqrt x}\)\)^{r+1}},
\tag4.5$$
where $U_n(x)$ denotes the $n$-th Chebyshev polynomial of the second
kind.
\endproclaim
\demo{Proof} By Lemma~$\Psi$, the permutations in the statement in
the theorem are in bijection with Dyck paths, which start at the
origin, return to the $x$-axis, and have exactly $r$ peaks at height
$k$. Such a Dyck path can be decomposed as
$$P_0u_1d_1P_1u_2d_2P_2\dots u_rd_rP_r,$$
where $P_0$ is the portion of $P$ from the origin until the first time
the height $k-1$ is reached, where the $u_i$'s are up-steps from
height $k-1$ to height $k$, where the $d_i$'s are down-steps from
height $k$ to height $k-1$, where the $P_i$'s,
$i=1,2,\dots,r-1$, are path portions which start and end at height
$k-1$ and never exceed height $k-1$, and where 
$P_r$ is the portion of $P$ from the last point at height $k-1$ until
the end of the path. Application of Theorem~A2 and use of Fact~A3
implies, by arguments that are more or less identically to those in
the proof of Theorem~3, that the generating function for those paths
is equal to
$$\frac {1} 
{\sqrt xU_{k}\(\frac {1} {2\sqrt x}\)}
\(\frac {U_{k-1}\(\frac {1} {2\sqrt x}\)} 
{\sqrt xU_{k}\(\frac {1} {2\sqrt x}\)}\)^{r-1}
\frac {1} 
{\sqrt xU_{k}\(\frac {1} {2\sqrt x}\)}x^{r},
$$
which simplifies to (4.5). \qed
\enddemo

The special case $k=3$ and $r=1$ of Theorem~10 appears, in an
equivalent form, in \cite{\RobeAA, Theorem~2}.

We could use an idea similar to the one in the proof of Theorem~3 to
express, for general $r$, 
the generating function for $123$-avoiding permutations 
with exactly $r$ occurrences of the pattern $(k-1)(k-2)\dots 1k$
in form of a sum, taken over all possible ways to arrange the peaks
that are at least at height $k$. However, it appears that it is not
possible to write the result in a way that is similarly elegant as
(3.5). However, for the asymptotics, the same reasoning as in the
proof of Theorem~5 remains valid. I.e., in this sum, the summand
which, asymptotically, provides the largest contribution, is again the
summand (4.5) (which is the same as the summand in (3.5) with
$\ell_1=r$ and all other $\ell_i$'s equal
to zero). Therefore an analogue of Theorem~5 in the present context
is true. More precisely, the number of $123$-avoiding permutations
with exactly $r$ 
occurrences of the pattern $(k-1)(k-2)\dots 1k$ is asymptotically as
large as the number of $132$-avoiding permutations with exactly $r$
occurrences of the pattern $12\dots k$. 
This extends \cite{\ChWeAA, Corollary~4.2}.

\proclaim{Theorem 11}
Let $r$ and $k$ be fixed nonnegative integers. 
Then, as $n$ becomes
large, the number of $123$-avoiding permutations with exactly $r$
occurrences of the pattern $(k-1)(k-2)\dots 1k$ is asymptotically
$$\(\frac {4\sin^2\frac {\pi} {k+1}} {k+1}\)^{r+1}\frac {n^r} {r!}
\(4\cos^2\frac {\pi} {k+1}\)^{n-r}\(1+O\(\frac {1} {n}\)\).\tag4.6$$ 
\line{\hfill \hbox{\qed}}
\endproclaim

\head Appendix. Generating functions for Motzkin and Dyck paths\endhead
A {\it Motzkin path} is a lattice path in the plane integer lattice $\Bbb
Z^2$, consisting of up-steps $(1,1)$, level-steps $(1,0)$, 
and down-steps $(1,-1)$, which never passes below the
$x$-axis. See Figure~2 for an example.

\vskip10pt
\vbox{
$$
\Gitter(12,5)(-1,0)
\Koordinatenachsen(12,5)(-1,0)
\Pfad(0,0),33141133443\endPfad
\DickPunkt(0,0)
\DickPunkt(11,2)
\hbox{\hskip6.5cm}
$$
\centerline{\eightpoint A Motzkin path}
\vskip7pt
\centerline{\eightpoint Figure 2}
}
\vskip10pt

Clearly, a Dyck path is just a Motzkin path without level-steps.

Given a Motzkin path $P$, we denote the {\it length\/} of the path
(i.e., the number of steps) by $\ell(P)$. Furthermore, 
we define the weight $w(P)$ of $P$ to be the
product of the weight of all its steps, where the weight of an
up-step is $1$ (hence, does not contribute anything to the weight), the
weight of a level-step at height $h$ is $b_h$, and the weight of
a down-step from height $h$ to $h-1$ is $\la_h$. Thus, the weight of
the Motzkin path in Figure~2 is
$b_2\la_2b_1b_1\la_3\la_2=b_1^2b_2\la_2^2\la_3$.

The theorem below, due to Flajolet, expresses the corresponding 
generating function
for all Motzkin paths which start at the origin and return to the
$x$-axis in form of a continued fraction.

\proclaim{Theorem A1}
{\smc(Flajolet \cite{\FlajAA, Theorem~1})}
With the weight $w$ defined as above, 
the generating function $\sum _{P} ^{}w(P)$, where the sum is over 
all Motzkin paths starting at the origin and returning to the
$x$-axis, is given by
$$
\cfrac 1\\1-b_0-
\cfrac \la_1\\1-b_1-
\cfrac \la_2\\1-b_2-\cdots\endcfrac\endcfrac\endcfrac\ .
\tag A.1$$
\line{\hfil \hbox{\qed}}
\endproclaim

Next we recall the expression, in terms of orthogonal
polynomials, for the generating function for Motzkin
paths in a strip. Although this is a result in the folklore of
combinatorics, probability, and statistics, the only explicit mention
that I am able to provide is \cite{\VienAE, Ch.~V, (27)}, which is a volume
that is not easily accessible. Therefore I include a sketch of proof.

\proclaim{Theorem A2} 
Define the sequence $(p_n(x))_{n\ge0}$ of polynomials
by
$$
xp_{n}(x)=p_{n+1}(x)+b_np_n(x)+\la_{n}p_{n-1}(x),\quad \quad 
\text{ for } n\geq 1,
\tag A.2$$
with initial conditions $p_0(x)=1$ and $p_1(x)=x-b_0$.
Furthermore, define\break $(Sp_n(x))_{n\ge0}$ to be the sequence of
polynomials which arises from the sequence $(p_n(x))$ by replacing $\la_i$ by
$\la_{i+1}$ and $b_i$ by $b_{i+1}$, $i=0,1,2,\dots$, everywhere in
the three-term recurrence {\rm(A.2)} and in the initial
conditions. Finally, given a polynomial $p(x)$ of degree $n$, we
denote the corresponding reciprocal polynomial
$x^np(1/x)$ by $p^*(x)$.

With the weight $w$ defined as before, the generating function $\sum
_{P} ^{}w(P)x^{\ell(P)}$, where the sum is over all
Motzkin paths which start at $(0,r)$, terminate at height $s$,
and do not pass above the line $y=K$, is given by
$$
\cases \dfrac {x^{s-r}p^*_r(x)S^{s+1}p_{K-s}^*(x)} {p_{K+1}^*(x)}
&r\le s,\\
\la_r\cdots\la_{s+1}\dfrac {x^{r-s}p^*_s(x)S^{r+1}p_{K-r}^*(x)}
{p_{K+1}^*(x)}
&r\ge s.
\endcases
\tag A.3
$$
\endproclaim
\demo{Sketch of Proof}
Motzkin paths which never exceed height $K$ 
correspond in a one-to-one fashion to 
walks on the path $P_{K+1}$ (this is the graph on the vertices
$v_0,v_1,\dots v_K$ where for $i=0,1,\dots, K-1$ the vertices $v_i$
and $v_{i+1}$ are connected by an edge, and there is a loop for each
vertex $v_i$). In this correspondence, an up-step from height $h$ to
$h+1$ in the Motzkin path corresponds to a step from vertex $v_h$ to
vertex $v_{h+1}$ in the walk, and similarly for level- and
down-steps. It is well-known (see e.g\. \cite{\StanAP, Theorem~4.7.2})
that the generating function for walks from $v_r$ to $v_s$ is given by 
$$\frac {(-1)^{r+s}\det(I-xA;r,s)} {\det(I-xA)},$$
where $A$ is the (weighted) {\it adjacency matrix} of $P_{K+1}$, where $I$
is the $(K+1)\times (K+1)$ identity matrix, and where $\det(I-xA;r,s)$
is the minor of $(I-xA)$ with the $r$-th row and $s$-th column
deleted. 

Now, the (weighted) adjacency matrix of $P_{K+1}$ with the property
that the weight of a particular walk would correspond to the weight $w$ of the
corresponding Motzkin path is the tridiagonal matrix
$$A=\pmatrix b_0&1&0&\dots\\
\la_1&b_1&1&0&\dots\\
0&\la_2&b_2&1&0&\dots\\
\vdots&\ddots&\ddots&\ddots&\ddots&\ddots&\vdots\\
&\dots&0&\la_{K-2}&b_{K-2}&1&0\\
&&\dots&0&\la_{K-1}&b_{K-1}&1\\
&&&\dots&0&\la_{K}&b_K\endpmatrix.$$
It is easily verified that, with this choice of $A$, we have
$\det(I-xA)=p_{K+1}^*(x)$ (by expanding the determinant 
with respect to the last row
and comparing with the three-term recurrence (A.2)), 
and, similarly, that the numerator in (A.3) agrees with  
$(-1)^{r+s}\det(I-xA;r,s)$.
\qed  
\enddemo

The special cases of Theorem~A2 in which $r=s=0$, respectively $r=0$
and $s=K$, appear also in \cite{\FlajAA, Sec.~3.1}.

The following is a well-known and easily verifiable fact:

\proclaim{Fact A3} If $b_i=0$ and $\la_i=1$ for all $i$, then the polynomials
$p_n(x)$ defined by the three-term recurrence {\rm(A.2)} are Chebyshev
polynomials of the second kind,
$$p_n(x)=U_n(x/2).\eqno\qed$$
\endproclaim

Fact~A3, in combination with Theorem~A2, tells that Chebyshev
polynomials of the second kind are intimately tied to the enumeration
of Dyck paths. 

Although we do not make use of it in this paper, we wish to emphasize 
that the enumeration of Motzkin 
paths (i.e., also allowing level-steps) is also intimately tied to
Chebyshev polynomials of the second kind.

\proclaim{Fact A4}
If $b_i=1$ and $\la_i=1$ for all $i$, then the polynomials
$p_n(x)$ defined by the three-term recurrence {\rm(A.2)} are also Chebyshev
polynomials of the second kind, namely
$$p_n(x)=U_n((x-1)/2).\eqno\qed$$
\endproclaim

Now we restrict our attention to Dyck paths. We refine the above
defined weight $w$ in the following way, 
so that in addition it also takes into account peaks: 
Given a Dyck path $P$, we define the weight $\hat w(P)$ of $P$ to be the
product of the weight of all its steps, where the weight of an
up-step is $1$, the weight of 
a down-step from height $h$ to $h-1$ which follows immediately after
an up-step (thus, together, form a peak of the
path) is $\nu_h$, and where the weight of
a down-step from height $h$ to $h-1$ which follows after another
down-step is $\la_h$. Thus, the weight of
the Dyck path in Figure~1 is 
$\nu_2\nu_4\nu_4\la_3\nu_3\la_2\la_1\nu_1=
\nu_1\nu_2\nu_3\nu_4^2\la_1\la_2\la_3.$

The theorem below, due to Roblet and Viennot, 
expresses the corresponding generating function
for all Dyck paths which start at the origin and return to the
$x$-axis in form of a continued fraction.

\proclaim{Theorem A5}
{\smc(Roblet and Viennot \cite{\RoViAA, Proposition~1})}
With the weight $\hat w$ defined as above, 
the generating function $\sum _{P} ^{}\hat w(P)$, where the sum is over 
all Dyck paths starting at the origin and returning to the
$x$-axis, is given by
$$
\cfrac 1\\1-(\nu_1-\la_1)-
\cfrac \la_1\\1-(\nu_2-\la_2)-
\cfrac \la_2\\1-(\nu_3-\la_3)-\cdots\endcfrac\endcfrac\endcfrac\ .
\tag A.4$$
\line{\hfil \hbox{\qed}}
\endproclaim

Finally, we remark that, from a different angle, 
Katzenbeisser and Panny \cite{\KaPaAD} have
undertaken an independent study of the enumeration of Motzkin paths.

\enddocument